\documentclass[reqno,11pt]{amsart}
\usepackage{amscd,amssymb,verbatim}
\usepackage{hyperref}\usepackage[small,nohug,heads=vee]{diagrams}
\diagramstyle[labelstyle=\scriptstyle]

\numberwithin{equation}{section}
\theoremstyle{plain}

\newcommand\alp{\alpha}         
\newcommand\bet{\beta}
         \newcommand\Gam{\Gamma}
         \newcommand\Del{\Delta}
\newcommand\eps{\varepsilon}

\newcommand\iot{\iota}

                \newcommand\Lam{\Lambda}
\newcommand\sig{\sigma}         \newcommand\Sig{\Sigma}

\newcommand\ome{\omega}         \newcommand\Ome{\Omega}

\newcommand\calD{{\mathcal{D}}}

\newcommand\calH{{\mathcal{H}}}

\newcommand\calO{{\mathcal{O}}}

\newcommand\bfv{{\mathbf v}}

\newcommand\RR{\mathbb{R}}

\newcommand\ZZ{\mathbb{Z}}

\newcommand\CC{\mathbb{C}}

\newcommand\NN{\mathbb{N}}

 \newcommand\grg{{\mathfrak{g}}}

\newcommand\nek{,\ldots,}
\newcommand\sdp{\times \hskip -0.3em {\raise 0.3ex
\hbox{$\scriptscriptstyle |$}}} 

\newcommand\End{\operatorname{End\,}}

\newcommand\Hom{\operatorname {Hom}}
\newcommand\Id{\operatorname {Id}}

\newcommand\IM{\operatorname{Im}}
\newcommand\Ind{\operatorname{Ind}}

\newcommand\Irr{\operatorname{Irr\, }}

\newcommand\Ker{\operatorname{Ker}}
\newcommand\Lie{\operatorname{Lie}}

\newcommand\supp{\operatorname{supp}}

\newcommand\hatR{{\widehat{R}}}

\renewcommand{\>}{\rangle}
\newcommand{\<}{\langle}

\theoremstyle{plain}
\newtheorem{Thm}[subsection]{Theorem}

\newtheorem{Prop}[subsection]{Proposition}

\newtheorem{Def}[subsection]{Definition}

\theoremstyle{remark}

\newtheorem{Rem}[subsection]{Remark}

\theoremstyle{remark}
\newtheorem{exam}[subsection]{Example}
\newtheorem{exercise}[subsection]{Exercise}

\usepackage{enumitem}
\usepackage{hyperref}
\usepackage[active]{srcltx}

\def\R{{ \! \rm \ I\!R}}

\newcommand{\Coker}{\operatorname{Coker}}
\newcommand{\Fred}{\operatorname{Fred}}
\newcommand{\tInd}{\operatorname{t-Ind}}
\newcommand{\n}{\nabla}
\begin{document}

\title{Index Theory of Non-compact $G$-manifolds}

\author{Maxim Braverman$^\dag$ and Leonardo Cano}
\address{Department of Mathematics\\
        Northeastern University   \\
        Boston, MA 02115 \\
        USA
         }
\address{Mathematics department, Andes University, Bogot\'a, Colombia.}
      
\thanks{${}^\dag$Supported in part by the NSF grant DMS-1005888.}


\maketitle


The index theorem, discovered by Atiyah and  Singer in
1963 \cite{AtSinger63}, is one of most important results in the
twentieth century mathematics. It found numerous applications in
analysis, geometry and physics. Since it was discovered numerous
attempts to generalize it were made, see for example
\cite{APS123,AtiyahTE,Atiyah76,Kasparov88,ConnesSkandalis84} to
mention a few; some of these generalizations gave rise to new very
productive areas of mathematics. In this lectures we first review
the classical Atiyah-Singer index theorem and  its generalization
to so called transversally elliptic operators~\cite{AtiyahTE} due
to Atiyah and Singer. Then we discuss the recent developments
aimed at generalization of the index theorem for transversally
elliptic operators to non-compact manifolds,
\cite{Paradan03,Br-index}.
\subsection*{Acknowledgments} This paper is  based on
the lectures given by the first author in  the Villa de Leyva
summer schools on  Geometric, Algebraic and Topological Methods for Quantum Field Theory.
 We are grateful to the students Camilo Orduz and Carlos Pinilla
   for careful note taking and helping us with preparing this manuscript.

 \section{The Fredholm Index}\label{S:Fredholm}
     
In this section we define the  index of an operator $A$ and discuss its main properties. The index is useful and nontrivial for operators defined on an  infinite dimensional vector space. But to explain the main idea of the definition let us start with the finite dimensional case.  

 \subsection{Finite dimensional case\label{SS:linearalgexample}}
 Consider a linear operator $A:\calH_{1}\rightarrow \calH_{2}$ between two {\em  finite dimensional} vector spaces, $\calH_{1}$ and $\calH_{2}$.  Then the {\em index} of $A$ is defined as 
 \begin{equation}\label{E:index}
   	\Ind(A)\ :=\ \dim \Ker(A)\ -\ \dim \Coker(A) \ \in  \ \ZZ.
 \end{equation}
where $\Coker(A):=\calH_{2}/\IM(A)$. Notice that, though $\dim\Ker(A)$ and $\dim\Coker(A)$ depend on
$A$, the index
\begin{equation}\label{E:indexfd}
   	\Ind(A)\ = \ \dim \calH_{1}\ - \ \dim \calH_{2}.
 \end{equation}
depends only on the spaces $\calH_1$ and $\calH_2$. 
    
\subsection{The Fredholm index}\label{SS:Fredholm}
Now suppose the spaces $\calH_1$ and $\calH_2$ are infinite dimensional Hilbert spaces. Then, in general, both $\Ker(A)$ and $\Coker(A)$ might be infinite--dimensional  and the index of $A$ cannot be defined. 

\begin{exercise} 
Show that if the dimension of $\Coker(A)$ is finite, then $\IM(A)$ is a closed subspace of $\calH_2$. 
\end{exercise}
  \begin{Def}\label{D:Fredholm}
Let $\calH_{1}, \calH_{2}$ be two Hilbert spaces and let $A:\calH_{1}\rightarrow \calH_{2}$ be  a bounded linear operator. 
We say that $A$ is a {\em Fredholm operator} iff $\dim\Ker(A)<\infty$, $\dim\Coker(A)<\infty$.
 \end{Def}

 If  $A:\calH_{1}\rightarrow \calH_{2}$ is a Fredholm operator, then its index $\Ind(A)$ can be defined by \eqref{E:index}. The notion of  index of a Fredholm operator was introduced by  Fritz Noether  \cite{Noether20},%
 \footnote{Fritz Noether was a fine mathematician with a very interesting and tragic biography. He was a  son of Max Noether and a younger brother of Emmy Noether. We refer the reader to \cite{Noether85,Noether90} for a short description of his biography.}%
 who also established the following ``stability" property of the index:

 \begin{Thm}\label{T:stability of index}
 \textbf{\em a.}\  If $A(t)$ is a Fredholm operator which depends continuously on the parameter $t$, then  $\Ind(A(t))$ is constant.
   \\
\textbf{\em b.}\  If $K$ is compact and $A$ Fredholm, then $\Ind(A+K)= \Ind(A)$.
\end{Thm}

   The following example illustrates the notion of index.

  \subsection{An example}\label{SS:Example shift}
 Let 
 \begin{equation}\label{E:l2}
	 \calH \ =\ l_2\ :=\ \big\{(x_1,x_2,\ldots\big|\, x_j\in \CC, \ \sum_{j=1}^\infty|x_j|^2<\infty\,\big\},
\end{equation}
and consider the linear operator $T:\calH\rightarrow \calH$  defined by
\begin{equation}\label{E:shiftby2}
  	T(x_{1},x_{2},x_{3},...) \ :=\ (x_{3},x_{4},...).
\end{equation}
 We have 
 \[
 	\Ker(T)\ =\ \big\lbrace (x_{1},x_{2},0,0,...) \in \calH\, \big|\,\, x_{1}, x_{2} \in  \mathbb{C}\,\big\rbrace.
\]
Thus $\dim\Ker(T)=2$. Since the image of $A$ is the whole space $\calH$, we conclude that  $\Coker(T)=\lbrace 0 \rbrace$. Thus $\Ind(T)=2$.
 
 \begin{Rem}\label{R:indexHHfd}
If  $\calH$ is a finite dimensional vector space then it follows from \eqref{E:indexfd} that $\Ind(A)=0$ for every linear operator  $A:\calH\to \calH$.  The example above shows that this is not true if $\dim{\calH}=\infty$. 
\end{Rem}

\begin{exercise}
Suppose $A:\calH_1\to \calH_2$ and $B:\calH_2\to \calH_3$ are Fredholm operators. Show that the operator $BA:\calH_1\to\calH_3$ is Fredholm and that 
\[
	\Ind (BA) \ = \ \Ind(A)\ +\ \Ind(B).
\]
\end{exercise}  
\begin{exercise}
Let $A:\calH_1\to \calH_2$ be a Fredholm operator. Show that the adjoint operator $A^*:\calH_2\to \calH_1$ is Fredholm and that 
\[
	\Ind(A)\ = \ -\ \Ind(A^*)\ = \ \dim\Ker(A)\ - \ \dim\Ker(A^*).
\]
\end{exercise}
\subsection{Application of index}\label{SS:aplindex}
One of the most common applications of index is based on Theorem~\ref{T:stability of index}. Suppose for example that $\Ind(A)>0$. This means that $\Ker(A)\not=\{0\}$, i.e., the equation $Ax=0$ has a non-trivial solution. Moreover, for any compact operator $K$ we have $\Ind(A+K)=\Ind(K)>0$. Hence, the equation 
\begin{equation}\label{E:A+Kx}
	(A+K)\,x\ =\ 0
\end{equation} 
also has a not-trivial solution. In applications to compute the index of an operator $B$ we often compute the index of a simpler operator $A$,  for which the kernel and cokernel can be explicitly computed, and then show that $B=A+K$ for $K$ a compact operator. If the index of $A$ is positive, we conclude that equation \eqref{E:A+Kx} has a non-trivial solution, even though we cannot find this solution explicitly.

\subsection{Connected components of the set of Fredholm operators}\label{SS:connectedFredholm}
Let $\Fred(\calH_1,\calH_2)$ denote the set of Fredholm operators $A:\calH_1\to \calH_2$. It is a metric space with the 
distance defined by $d(A,B):=\|A-B\|$. Theorem~\ref{T:stability of index} implies that if operators $A$ and $B$ belong 
to the same connected component of $\Fred(\calH_1,\calH_2)$, then $\Ind(A)= \Ind(B)$.  
In other words, $\Ind(A)$ is an invariant of the connected component of $\Fred(\calH_1,\calH_2)$.
In fact, $\Ind(A)$ determines the connected component of $A$ in
$\Fred(\calH_1,\calH_2)$ (see \cite[proposition~7.1]{Lawson}).

\subsection{The group action}\label{SS:group action}
The index $\Ind(A)$ of a Fredholm operator $A$ is an integer. If a compact group $G$ acts on a Hilbert space $\calH$ 
one can define a richer invariant, as we shall now explain. 

Recall that a $G$-representation $V$ is  called {\em irreducible} if it has no non-trivial  $G$-invariant subspaces. 
We denote by $\Irr(G)$ the set of irreducible representations of $G$. The next theorem, cf. \cite{BroeckerTammo}, shows that they are the 
building-blocks  of all the other representations.

\begin{Thm}\label{T:irrdecomp}
Any finite dimensional representation $U$ of $G$ has a unique decomposition into a sum of irreducible representations:
\begin{equation*}
	U\ = \ \bigoplus_{V\in\Irr(G)}m_VV.
\end{equation*}
Here $m_VV$ stands for the direct sum $V\oplus V\oplus\cdots\oplus V$ of $m_V$copies of $V$.
 \end{Thm}
 The numbers $m_V\in \NN$ are called the {\em multiplicities of the irreducible representation $V$ in $U$}. 
\begin{Def}\label{D:equivariant}
 Let $G$ be a compact group acting on the vector spaces $\calH_{1}$ and $\calH_{2}$. A linear transformation $A:\calH_{1}\rightarrow \calH_{2}$ is $G$-equivariant  iff {\em it commutes with the action of the group}, i.e.,
 \[
 	gAx\ =\ Agx, \qquad \text{for all}\quad x\in \calH_1. 
\]
\end{Def}

\begin{exercise}\label{Exer:equivariant kernel}
Let $G$ be a compact group acting on Hilbert spaces $\calH_1$ and $\calH_2$. Let $A:\calH_1\to \calH_2$ be a $G$-equivariant  Fredholm operator.

a.\ Show that $\Ker(A)$ is invariant under the action of $G$, i.e. for  any $x\in \Ker(A), \ g\in G$ we have $gx\in \Ker(A)$.
 Thus the restriction of the action of $G$ on $\calH_1$ to $\Ker(A)$ defines a $G$-action on $\Ker(A)$.

b.\  Define a natural action of $G$ on $\Coker(A)$. 
\end{exercise}

 The above exercise and Theorem~\ref{T:irrdecomp} show that there is a unique decomposition
 \[ 
   \Ker(A)\ =\ \bigoplus_{V\in \Irr(G)} m^{+}_V\,V;\qquad
    \Coker(A)\ =\ \bigoplus_{V\in \Irr(G)} m^{-}_V\,V.
 \]
 
 \begin{Rem}\label{R:finitesum}
 Note the though the set $\Irr(G)$ is infinite both sums above are actually finite, since only finitely many numbers $m_+,m_-$ can be non-zero.
 \end{Rem}
\subsection{The ring of characters}\label{SS:characters}
Under the direct sum, $\oplus$, the set of finite dimensional representations of the group $G$ form an abelian semigroup with
identity. Thus we can associate to it the {\em Grothendieck group $R(G)$} by considering formal differences $V-U$ of finite
dimensional representations of $G$. Let us give a formal definition:

\begin{Def}\label{D:R(G)}
Consider the set of formal differences, $V_{1}-V_{2}$, of finite dimensional representations $V_1$ and $V_2$ of $G$. The quotient of this set by the equivalence relation\footnote{In general, the definition of the equivalence relation in the Grothendieck group is slightly more complicated. But, in our case, it is equivalent to the one given here.}
\[
	V-U\ \sim\  A-B\ \Longleftrightarrow\ V \oplus B \simeq A\oplus U,
\]
 is an abelian group denoted by $R(G)$ and called the {\em group of characters of $G$}. Formally, 
 \[
 	R(G)\ := \ \big\{(V,U)\big|\, V,U\ \text{\em are finite dimensional representations of}\ G\,\big\}/\sim.
\]
\end{Def}

The tensor product of representations defines a product on $R(G)$ and thus defines a ring structure on $R(G)$. That is why $R(G)$ is called the {\em ring of characters of $G$}. We will not use the ring structure on $R(G)$ in this lectures and will not give a formal definition. 

Using Theorem~\ref{T:irrdecomp} one easily obtains the following alternative description of $R(G)$:
\begin{Prop}\label{P:R(G)=sum}
As an abelian group the ring  of characters $R(G)$ is isomorphic to the free abelian group generated by the set $\Irr(G)$ of irreducible representations of $G$, i.e., to the group of expressions
\begin{equation}\label{E:R(G)=sum}
	\Big\{\bigoplus_{V\in \Irr(G)}\,m_VV\big|\, m_V\in \ZZ, \ \text{\em only finitely many numbers}\ m_V\not=0\,\Big\}.
\end{equation}
\end{Prop}

\begin{exercise}\label{Exer:ringstructureonR(G)}
Introduce a ring structure on the group \eqref{E:R(G)=sum} and show that the isomorphism of Proposition~\ref{P:R(G)=sum} is an isomorphism of rings. 
\end{exercise}
\begin{Rem}\label{R:hatR}
If one  allows infinite formal sums in \eqref{E:R(G)=sum} then one obtains the definition  of the  {\em completed ring of characters} of $G$:
\begin{equation}\label{E:hatR(G)=sum}
	\hatR(G)\ = \ \Big\{\bigoplus_{V\in \Irr(G)}\,m_VV\big|\, m_V\in \ZZ, \,\Big\}.
\end{equation}
This ring will play an important role in Sections~\ref{S:trelliptic} and \ref{S:indexnc}.
\end{Rem}

\subsection{The equivariant index}\label{SS:eqindex}
Suppose a compact group  $G$ acts  on Hilbert spaces $\calH_{1}$ and $\calH_{2}$ and let $A:\calH_{1}\rightarrow \calH_{2}$ be a $G$-equivariant Fredholm operator.

\begin{Def}\label{D:eqindex}
The {\em equivariant index} $\Ind_G(A)$ of $A$ is defined to be the formal difference between $\Ker(A)$ and $\Coker(A)$ in $R(G)$:
\begin{equation}\label{E:eqindex}
	\Ind_G(A)\ := \ \Ker(A)\ -\ \Coker(A)\ = \ \bigoplus_{V\in\Irr(G)}(m^+_V-m^-_V)\,V\ \in \ R(G).
\end{equation}
\end{Def}
By Remark~\ref{R:finitesum} the sum in the right hand side of \eqref{E:eqindex} is finite.

\begin{exam}\label{Exam:eqequation}
Let $V_{0}$ denote the trivial representation. If the inequality $m^{+}_{V_0}-m^{-}_{V_0} > 0$ holds, that means that there  exist a nontrivial  solution of the equation $Ax=0$ which is invariant under the action of $G$.
\end{exam}

\begin{exercise} 
Let 
\[
	G\ :=\ \mathbb{Z}_{2}=\lbrace 1;-1\rbrace.
\]
To avoid a confusion we denote by $q$  the non-trivial element of $G$ (thus $q=-1$). We define an action of $G$ 
on $l^{2}$ by
 \[
 	q(x_{1},x_{2},...)\ =\ (x_{2},x_{1},x_{4},x_{3},...).
\]
Observe that $G$ has two irreducible representations, namely: $V_{0}=\mathbb{C}$ with the trivial action, and $V_{1}=\mathbb{C}$ with the action $qz=-z$. 

Show that the operator \eqref{E:shiftby2} is $G$-equivariant and that  $\Ker(T)=V_{0}\oplus V_{1}$. Use this result to compute $\Ind_G(T)$. 
 \end{exercise}

 \section{Differential operators}\label{S:differential}
  
 In this section we introduce differential operators on manifolds and discuss their main properties. 
We also define elliptic differential operators. If the manifold is compact, then any elliptic operator on it is Fredholm. 
In the next section we will discuss the index of elliptic 
differential operators and the Atiyah-Singer theorem, which computes this index in terms of topological 
information of the manifold.

  \subsection{Differential operators\label{SS:difoper}}
 We recall the definition of linear differential operators in $\mathbb{R}^{n}$. Let $D_j:C^\infty(\RR^n)\to C^\infty(\RR^n)$ denote the operator 
 \[
 	D_jf\ = \ \frac1i\,\frac{\partial f}{\partial x_j}, \qquad j=1\nek n.
\]
 
 A {\em multiindex} is an $n$-tuple of non-negative integers $\alp=(\alp_1\nek \alp_n)\in \ZZ^n_+$. For a multiindex $\alp$ we set 
 \[
 	|\alp|\ := \ \alp_1+\cdots+\alp_n
\]
and
\[
	D^\alp \ := \ D_1^{\alp_1}D_2^{\alp_1}\cdots D_n^{\alp_n}:\, C^\infty(\RR^n) \ \longrightarrow \ C^\infty(\RR^n).
\]
 
 \begin{Def}\label{D:Difoper}
 A {\em linear differential operator of order $k$}, is a linear operator $\calD$ of the form
 \begin{equation}\label{E:diffoperRn}
 	 \calD\ =\   \sum_{|\alpha|\leq k} a_{\alpha}(x) \, D^\alp:\, C^\infty(\RR^n) \ \longrightarrow \ C^\infty(\RR^n) ,
\end{equation}
where  $a_{\alpha}(x)\in C^{\infty}(\mathbb{R}^{n})$.
 \end{Def}

Let us denote $\<\cdot, \cdot\>$ the scalar product on $C^\infty_c(\R^n)$ defined by
\begin{equation*}
\<f, g\>:=\int_{\R^n} f(x)\overline{g}(x)dx.
\end{equation*}
\begin{Def}\label{D:Sobolev space}
Let $C^\infty_c(\RR^n)$ denote the set of smooth functions with compact support on $\RR^n$. For $k \in \NN$ define the scalar product on $C^\infty_c(\RR^n)$ by the formula
\begin{equation}\label{E:Sobolevproduct}
	\<f,g\>_k\ := \ \sum_{|\alp|=k}\, \<D^\alp f,D^\alp g\>.
\end{equation}
The corresponding norm 
\begin{equation}\label{E:Sobolevnorm}
	\|f\|_k\ := \ \sqrt{\<f,f\>_k}
\end{equation}
is called the {\em $k$-th Sobolev norm}. 

The {\em Sobolev space $H^{k}(\RR^n)\subset L^2(\RR^n)$} is the completion of the space $C^\infty(\RR^n)$ with respect to the norm \eqref{E:Sobolevnorm}.
 \end{Def}
 The following theorem, \cite[Corollary~3.19]{Adams-SobolevSpaces}, gives an alternative description of the Sobolev space $H^{k}(\RR^n)$:
 
 \begin{Thm}\label{T:Sobolevcompl}
The Sobolev space $H^{k}(\RR^n)$ is equal to the space of square integrable functions whose distributional derivatives up to order $k$ are in $L^{2}(\mathbb{R}^{n})$. More explicitly
\[ 
  	H^{k}(\RR^n)\\ :=\ 
	\big\{ f: \RR^{n}\to\CC\,\big|\, D^\alp f\in L^2(\RR^n) \ \text{\em for all}  \
		\alp\in\ZZ_+^n\ \text{\em with}\  |\alp|\le k\, \big\}.
\]
 \end{Thm}

\begin{Rem}
If $\Ome\subset \RR^n$ is an open set one can define the Sobolev space $H^k(\Ome)$ of functions on $\Ome$. An analogue of Theorem~\ref{T:Sobolevcompl} still holds if the boundary of $\Ome$ is sufficiently nice, see \cite[Theorme~3.18]{Adams-SobolevSpaces}.
\end{Rem}

\begin{exercise}\label{Exer:Extension to Sobolev}
Let $K\subset M$ be a compact set. Denote by 
\[
	H^k_K(\RR^n)\ :=\ \big\{\,f\in H^k(\RR^n)\,\big|\,\supp(f)\subset K\,\big\}.
\]
For every $m\ge k$ a differential operator $\calD:C^\infty_c(\R^n)\to C^\infty_c(\R^n)$ of order $k$ can be extended continuously to a bounded operator $\calD:H^{m}_K(\mathbb{R}^{n})\rightarrow L^{2}(\mathbb{R}^{n})$.
\end{exercise}
With a bit more work one can show the following
\begin{Thm}{\upshape \cite[Theorem~6.2]{Adams-SobolevSpaces}}\label{T:Rellich}\textbf{\em (Rellich's lemma)}
If $K\subset \RR^n$ is a compact set, then the embedding
$H^m_K(\RR^n)\hookrightarrow H^k(\RR^n)$ is a compact operator for
any $m>k$.
\end{Thm}
\begin{exercise}\label{Exer:compact diff oper}
Use Rellich's lemma to show that  the restriction of a differential operator $\calD$ of order $k<m$  to $H^m_K(\RR^n)$ defines a compact operator $H^m_K(\RR^n)\to L^2(\RR^n)$.
\end{exercise}

\subsection{Matrix-valued differential operators}\label{SS:matrixvalued}
We extend Definition ~\ref{D:Difoper} to operators acting on vector valued functions on $\RR^n$. Let $C^\infty(\RR^n,\RR^N)$ denote the set of smooth vector valued functions 
\[
	f\ = \ (f_1\nek f_N):\,\RR^n\ \longrightarrow \ \RR^N, \qquad f_j\in C^\infty(\RR^n).
\]
Using the scalar product of $\RR^N$ we can generalize the constructions of the previous section to define the spaces 
$L^2(\RR^n,\RR^N), H^k(\RR^n,\RR^N)$ and $H_K^k(\RR^n,\RR^N)$, etc.

An operator $\calD:\,C^\infty(\RR^n,\RR^{N_1})\to C^\infty(\RR^n,\RR^{N_2})$ is a (matrix-valued) differential operator of order $k$ if it is given by the formula \eqref{E:diffoperRn} with 
\[
	a_\alp(x)\ =\ \big\{a_{\alp,ij}(x)\big\}_{\tiny\substack{1\le i\le N_2\\ 1\le j\le N_1}}\in \operatorname{Mat}_{N_1\times N_2},
	\qquad a_{\alp,ij}(x)\in C^\infty(\RR^n).
\]
As above, a matrix-valued differential operator $\calD$ of order $k$ defines a bounded operator 
\[
	\calD:\,H^k(\RR^n,\RR^N)\ \longrightarrow\  L^2(\RR^n,\RR^N)
\]
and a compact operator
\[
	\calD:\,H^m_K(\RR^n,\RR^N)\ \longrightarrow\  L^2(\RR^n,\RR^N), \qquad \text{for}\ \ m>k.
\]

\subsection{Vector bundles}\label{SS:bundles}
Consider a smooth manifold $M$and let $E$ be a vector bundle over $M$. We refer the reader to \cite{Taubes-DifGeom} for the definition and the basic properties of vector bundles. To understand what follows one should remember that $E$ is itself a smooth manifold which is equipped with a {\em projection map} $\pi:E\to M$. For each $x\in M$ the preimage $\pi^{-1}(x)$ is called the {\em fiber of $E$ over $x$} and is denoted by $E_x$. It is assumed that each fiber $E_x$ is isomorphic to the linear space $\CC^N$ (however no preferred isomorphism is fixed in general). 

\begin{exam}
a. The {\em trivial bundle} $E=M\times\CC^N$. Here $\pi:E\to M$ is the projection on the first factor. In this case each fiber $E_x$  is canonically isomorphic with $\CC^N$.

b. The {\em tangent bundle} $TM$. This is a vector bundle whose fiber at each point $x\in M$ is equal to the tangent space $T_xM$. 
\end{exam}

A {\em smooth section} of a vector bundle $E$ is a smooth map $s:M\to E$ such that 
\[
	\pi\big(s(x)\big) \ = \  x, \qquad \text{for all}\ \ x\in M. 
\] 
In other words for each $x\in M$ we define an element $s(x)\in E_x$. Notice that for the case of the trivial bundle $E=M\times\CC^N$ defining a section $s$ is the same a defining a vectorÐvalued function $M\to \CC^N$. Thus the notion of smooth sections generalizes the notion of smooth vectorÐvalued functions on $M$. 

We denote the set of smooth sections of $E$ by $C^\infty(M,E)$.

\subsection{Differential operators on manifolds}\label{SS:operators on manifolds}
Let $M$ be a smooth manifold and let $E$ and $F$  be  vector bundles over $M$. We say that a linear operator 
\begin{equation}\label{E:DEtoF}
	\calD:\,C^\infty(M,E)\ \longrightarrow\ C^\infty(M,F)
\end{equation}
is {\em local} if for every section $f\in C^\infty(M,E)$, we have $\supp(\calD{}f)\subset \supp(f)$. In this case for 
any open cover $\{U_i\}_{i=1}^m$ of $M$ the operator $\calD$ is completely determined by its restriction to functions 
with supports in one of $U_i$. If all the sets $U_i$ are sufficiently small, we can fix coordinates 
\[
	\phi_i:\,U_i\ \overset{\sim}\longrightarrow\ V_i\subset \RR^n \qquad (i=1\nek m)
\] 
in $U_i$ and also trivializations of $E$ and $F$ over $U_i$. Then the restriction of $\calD$ to sections supported 
in $U_i$ can be identified with a map
\[
	\calD_i:\, C^\infty(V_i,\RR^{N_1})\ \longrightarrow\ C^\infty(V_i,\RR^{N_2}).
\]

\begin{Def}\label{D:diffoperM}
A local linear operator $\calD:C^\infty(M,E)\to C^\infty(M,F)$ is called a differential operator of order $k$ if  one can choose a cover $\{U_i\}_{i=1}^m$ of $M$  and coordinate systems on $U_i$ such that for every $i=1\nek m$ the operator $\calD_i$ has the form \eqref{E:diffoperRn}.
\end{Def}

\subsection{Sobolev spaces of sections}\label{SS:Sobolevsection}
We shall now introduce the Sobolev spaces of sections of a vector bundle $E$.  Let us choose a measure on $M$ and a 
scalar product on the fibers of $E$. Then we can consider the space $L^2(M,E)$ of square-integrable section of $E$. 
It is not hard to show that if the manifold $M$ is closed (i.e
compact and without boundary) the space $L^2(M,E)$ is independent of the choice of the 
measure and of  the scalar product on $E$. 

One can define the spaces $H^k(M,E)$ of Sobolev sections of $E$ in a way similar to Definition~\ref{D:Sobolev space}. Roughly one users a partition of unity to define a scalar product $\<\cdot,\cdot\>_k$ on $C^\infty_k(M,E)$ as a combination of scalar products \eqref{E:Sobolevproduct} on each coordinate neighborhood. Then the Sobolev space $H^k(M,E)$ is the completion of  $C^\infty_c(M,E)$ with respect to the norm defined by this scalar product. We refer to \cite{Taubes-DifGeom} for  details. As in the case $M=\RR^n$ we have the following:
\begin{Thm}\label{T:DHtoL}
Let  $\calD:C^\infty(M,E)\to C^\infty(M,F)$ be a differential operator of order $k$.  

\textbf{a.} \ For every $m\ge k$ and every compact $K\subset M$, the operator $\calD$ extends continuously to an operator 
\[
	\calD:\,H^m_K(M,E)\ \longrightarrow\ L^2(M,F).
\] 

\textbf{b.}\ If $m>k$ then for every compact $K\subset M$ the operator $\calD:H^m_K(M,E)\to L^2(M,F)$ is compact. 
\end{Thm}

\subsection{The symbol of a differential operator}\label{SS:symbol}
We now define the notion of a leading symbol of a differential operator which is crucial for the discussion of the Atiyah-Singer index theorem. 

\begin{exercise}\label{Exer:leadingsymbol}
{a.}\ Let $\calD:C^\infty(M,E)\to C^\infty(M,F)$ be a differential operator of order $k$.  Fix $x_0\in M$, a cotangent vector $\xi\in T_{x_0}^*M$ and a vector $e\in E_{x_0}$. Let $f\in C^\infty(M)$ be a smooth function, such that $df_{x_0}=\xi$ and let $s\in C^\infty(M,E)$ be a section, such that $s(x_0)=e$. Set 
\begin{equation}\label{E:leadingsymbol}
	\sig_L(\calD)(x_0,\xi)\,e \ := \ \lim_{t\to \infty}\, t^{-k}\,\calD\big(e^{itf(x)}s(x)\big)\big|_{x=x_0}.
\end{equation}
Show that \eqref{E:leadingsymbol} is independent of the choice of $f$ and $s$ and that $\sig_L(\calD)(x_0,\xi)\,e$ is linear in $e$. Thus we can view $\sig_L(\calD)(x_0,\xi)$ as a linear map $E_{x_0}\to F_{x_0}$. 

{b.}\ Suppose that $\calD$ is given in local coordinates by the formula \eqref{E:diffoperRn}. Using local coordinates we can identify $T_{x_0}M$ with $\RR^n$. Then 
\begin{equation}\label{E:symbolincoord}
	\sig_L(\calD)(x_0,\xi)\ = \ \sum_{|\alp|=k}\,a_\alp(x_0)\,\xi_1^{\alp_1}\cdots\xi_n^{\alp_n}.
\end{equation}
\end{exercise}
Notice, that the above exercise implies that as a function on $T^*M$ the right hand side of \eqref{E:symbolincoord} is independent of the choice of a coordinate system near $x_0$. 

\begin{Def}\label{D:leadingsymbol}
Let $\calD:C^\infty(M,E)\to C^\infty(M,F)$ be a differential operator of order $k$. The {\em leading symbol} $\sig_L(\calD)(x_0,\xi)$ ($x_0\in M,\ \xi\in T^*_{x_0}M$) of $\calD$ is an element of  $\Hom(E_x,F_x)$ defined by \eqref{E:leadingsymbol}.
\end{Def}

\begin{Rem}
Using the local coordinate representation \eqref{E:diffoperRn} of $\calD$ on can also define the {\em full symbol} 
\[
	\sig(\calD)(x_0,\xi)\ = \ \sum_{|\alp|\le k}\,a_\alp(x_0)\,\xi_1^{\alp_1}\cdots\xi_n^{\alp_n}.
\]
But as opposed to the leading symbol the full symbol depends on the choice of the local coordinates.
\end{Rem}

\subsection{The leading symbol as a section of the pull-back bundle}\label{SS:symbolsection}
Let $\pi:T^*M\to M$ denote the projection. For a vector bundle $E$ we denote by $\pi^*E$ the {\em pull-back} bundle over $T^*M$. This is a vector bundle over $T^*M$ whose fiber over $(x,\xi)\in T^*M$ is isomorphic to the fiber $E_x$ of $E$, see \cite[\S5.1]{Taubes-DifGeom} for a precise definition of the pull-back bundle. 

Let $\calD:C^\infty(M,E)\to C^\infty(M,F)$ be a differential operator. Then for $(x,\xi)\in T^*M$ , the leading symbol 
$\sig_L(\calD)(x,\xi)$ is a linear map 
\[
	\sig_L(\calD)(x,\xi):\,\pi^*E_{(x,\xi)}\ \to\ \pi^*F_{(x,\xi)}.
\] 
Hence, we can view $\sig_L(\calD)$ as a section of the vector bundle $\Hom(\pi^*E,\pi^*F)$.

\subsection{Elliptic differential operators}\label{SS:elliptic}

\begin{Def}\label{D:elliptic}
A differential operator $\mathcal{D}:C^\infty(M,E) \to C^\infty(M,F)$ is called {\em elliptic} if $\sig_L(\calD)(x,\xi)$ is invertible for all $\xi\not=0$. 
\end{Def}

\begin{exam}
Consider the  {\em Laplace operator} 
\[
	\Delta\ =\ - \frac{\partial^2}{\partial x_1^2}- \frac{\partial^2}{\partial x_2^2}-\cdots-\frac{\partial^2}{\partial x_n^2}
\] 
on $\RR^n$.  Its leading symbol $\sigma_L(\Delta)$ is given by the formula
\[
    \sigma_L(\Delta)(\xi_1,\xi_2,\nek \xi_n) \ = \ \xi_1^2+\xi_2^2+\cdots+\xi_n^2
\]
is invertible for $\xi=(\xi_1\nek\xi_n)\not=0$. Hence, $\Del$ is elliptic. 
\end{exam}

 Elliptic operators play a very important role in analysis and in the theory of index. The following results, cf. \cite{ShubinPDObook}, describes some of the main properties of elliptic operators

\begin{Thm}\label{T:elliptic}
Suppose  $\mathcal{D}:C^\infty(M,E) \to C^\infty(M,F)$ is an elliptic differential operator. Then
\\
\textbf{a.}\ (Elliptic regularity) If\, $\mathcal{D}f=u$ and $u \in
C^\infty(M,F)$ then  $f \in C^\infty(M,E)$.
\\
\textbf{b.}\ If the manifold  $M$ is compact, then the operator $\mathcal{D}$ is a Fredholm.
\end{Thm}


Part b. of  the above theorem implies that the index $\Ind(\calD)$ of an elliptic operator $\calD$ on a compact manifold is defined, cf. Section~\ref{S:Fredholm}.  In the next section we will discuss the Atiyah-Singer index theorem, which computes this index. We finish this section with the following
\begin{exercise}\label{Exer:indexlsymbol}
Use Theorem~\ref{T:DHtoL} to show that $\Ind(\calD)$ depends only on the leading symbol of $\calD$, i.e. if $\sig_L(\calD_1)= \sig_L(\calD_2)$ then $\Ind(\calD_1)= \Ind(\calD_2)$.
\end{exercise}

\begin{exercise}
Let $M=S^1$ be the circle and consider the operator
\[
	\calD\ := \ -i\,\frac{d}{dt}+\sin t:\, C^\infty(M)\ \to \ C^\infty(M).
\]
Compute $\Ind(\calD)$.

\noindent{\em Hint}:\  Compute the index of a simpler operator with the same leading symbol. 
\end{exercise}
\section{The Atiyah-Singer index theorem}\label{S:Atiyah-Singer}

\newcommand{\SEll}{\operatorname{SEll}}

In this section we present a $K$-theoretical formulation of the Atiyah-Singer index theorem.

\subsection{Index as a topological invariant}\label{SS:homotinvindex}
At the end of the last section (cf. Exercise~\ref{Exer:indexlsymbol}) we saw that the index of an elliptic differential operator on a compact manifold depends only on its leading symbol. Moreover, Theorem~\ref{T:stability of index}.a implies that the index does not change when we deform the leading symbol. More precisely, let us consider the space $\SEll(E,F)$ of smooth sections of $\sig(x,\xi)\in \Hom(\pi^*E,\pi^*F)$, which are invertible for $\xi\not=0$. We refer to  $\SEll(E,F)$ as the {\em space of elliptic symbols}. We endow it with the topology of uniform convergence on compact sets. 

\begin{Thm}\label{T:homotinvindex}
Let $E$ and $F$ be vector bundles over a compact manifold $M$. The index of an elliptic differential operator 
$\calD:C^{\infty}(M,E)\to C^{\infty}(M,F)$ is determined by the connected component of\, $\SEll(E,F)$ in which the leading symbol $\sig_L(\calD)$ of $\calD$ lies. 
\end{Thm}

This result shows that $\Ind(\calD)$ is a topological invariant.  About 50 years ago Israel Gel'fand 
\cite{Gelfand57Congress,Gelfand60} formulated a problem: how to compute the index of an elliptic operator using 
only its
leading symbol. This problem was solved brilliantly by Atiyah and Singer \cite{AtSinger63,AtSinger1} 
(see also, for example, \cite{Palais65,BeGeVe}). They did more than just a calculation of the index. They associated 
to each elliptic symbol an element of so called $K$-theory (see Section~\ref{SS:Ktheory} below). Then
they associated a number -- {\em  the topological index $\tInd(\calD)$} -- to each element of the $K$-theory. Schematically, their construction can be expressed as
\begin{equation}\label{E:ASscheme}
   \mathcal{D} \ \rightsquigarrow \ \sigma_L(\mathcal{D}) \ \rightsquigarrow \
    \text{an element of $K$-theory} \ \rightsquigarrow \ \tInd(\sig_L(\calD)))\ \in \ \ZZ.
\end{equation}
The composition of the arrows in the above diagram leads to a map
\[
    \text{elliptic operators}  \ \rightarrow  \ \ZZ, \qquad \ \mathcal{D} \mapsto \ \tInd(\sigma_L(\mathcal{D})).
\]
which is called the topological index. Note, that, as we will explain below,  the topological index is constructed using purely topological methods, without any analysis involved.

The following result is the simplest form of the Atiyah-Singer index theorem

\begin{Thm}\label{T:index}\textbf{\em (Atiyah-Singer)}
$\Ind (\calD) \ =\  \tInd(\sigma_L (\mathcal{D}))$.
\end{Thm}

\begin{Rem}
Though from Exercise~\ref{Exer:indexlsymbol} we know that the index of $\mathcal{D}$ could be computed out of the leading 
symbol $\sigma_L (\mathcal{D})$, we call the reader to appreciate the power of the above theorem. The index of 
$\mathcal{D}$ gives you an information about the
kernel and the cokernel of $\mathcal{D}$, i.e., about the spaces of solutions of differential equations 
$\mathcal{D}f=0$ and $\mathcal{D}^*u=0$. The index theorem allows to obtain this information 
{\em without solving the differential equations} by purely topological methods.
\end{Rem}

We shall now explain the meaning of the terms in \eqref{E:ASscheme}. First, we will briefly review the notion of the (topological)
$K$-theory.
\subsection{$K$-theory}\label{SS:Ktheory}
The direct sum $\oplus$ defines a structure of a semi-group on the  set of vector bundles 
ovex $X$. If $X$ is a compact manifold. then the $K$-theory of $X$ is just the Grothendieck group of this semi-group. 
In general the definition is a little bit more complicated. This is, essentially, because we are interested in the 
$K$-theory {\em with compact support}. In fact, one of the definitions of the $K$-theory of a non-compact manifold $X$ 
is the $K$-theory of the {\em one-point compactification of $X$}. It will be more convenient for us to use an equivalent 
definition, which is based on considering bundle maps
\begin{center}
\setlength{\unitlength}{.2in}
\begin{picture}(8,6)
\put(0,4){$E$} \put(4.5,4){$F$}  \put(2.1,0){$X$} \put(1,4.3){\vector(1,0){3.2}} \put(0.5,3.7){\vector(1,-2){1.4}} \put(4.5,3.7){\vector(-1,-2){1.4}}
\put(2,5){$\sigma$}
\end{picture}
\end{center}
such that the induced map of fibers  $\sigma(x):E_x\to F_x$ is invertible for all $x$ outside of a compact set $K\subset X$. Two such maps $\sigma_1:E_1\to F_1$
and $\sigma_2:E_2\to F_2$ are said to be equivalent if there exist integers $k_1,k_2\ge 0$ such that the maps
\begin{equation}\label{E:Kequivalence}
    \sigma_1\oplus \Id:\, E_1\oplus \CC^{k_1}\ \to\ F_1\oplus \CC^{k_1}, \quad\text{and}\quad
    \sigma_2\oplus\Id:\,E_2\oplus \CC^{k_2}\ \to \ F_2\oplus \CC^{k_2},
\end{equation}
are {\em homotopic in the class of maps invertible outside of a compact set}. In this case we write $\sigma_1\thicksim \sigma_2$.

\begin{Def}\label{D:Ktheory}
Let $X$ be a topological space. The \emph{K-theory} $K(X)$ of $X$ is defined by
\[
    K(X) \ = \ \{\sigma: E \rightarrow F \mid \sigma(x) \textrm{ is invertible outside of a compact set} \} / \thicksim
\]
\end{Def}

\begin{Rem}
If $X$ is compact then any two maps $\sigma_1,\sigma_2:E\to F$ between the same vector bundles are equivalent. Thus $K(X)$ can be described as the set of
pairs of vector bundles $(E,F)$ subject to an appropriate equivalence relation. We think about this pair as about {\em formal difference} of the
bundles $E$ and $F$, and we, usually, denote this pair by $E-F$.
\end{Rem}

The direct sum of vector bundles defines a structure of an abelian group on $K(X)$. The tensor product defines a multiplication. Together these two operations make $K(X)$ a ring. 

\begin{exam}
 If $X=\{pt\}$ then a vector bundle over $M$ is just a vector space and an element of $K(X)$ is  a pair of two vector
spaces $(E,F)$ up to an equivalence. In this case one easily sees from \eqref{E:Kequivalence} that  the only invariant of the pair $(E,F)$ is the number
\[
    \dim E \ - \ \dim F \ \in \ \ZZ.
\]
Thus
\begin{equation*}
    K(\{pt\})\ \cong\ \mathbb{Z}.
\end{equation*}
\end{exam}

\begin{exam} \label{Exam:siginK}
Let $\calD:C^\infty(M,E)\to C^\infty(M,F)$ be an elliptic differential operator. Then for each $x\in M$ and each $\xi\in T^*_xM$ we have $\sigma_L (\mathcal{D})(x,\xi): E_x\to F_x$. In other words, $\sigma(\mathcal{D})$ defines a bundle map
\begin{equation}\label{E:sig-Kth}
    \sigma_L (\mathcal{D}): \pi^\star E~~\rightarrow~~\pi ^\star F.
\end{equation}
If $M$ is a {\em compact manifold}, then the ellipticity condition implies that the map \eqref{E:sig-Kth} is invertible outside of a compact subset of $X=T^*M$. Hence, \eqref{E:sig-Kth}  defines an element of $K(T^*M)$.

This example shows the relevance of the K-theory to the problem of index. 
\end{exam}

\subsection{The push-forward map in K-theory}\label{SS:pushforwardK}

One of the most important facts about the $K$-theory is that it has all the properties of a cohomology theory with 
compact supports. In particular, given
an inclusion $j:Y\hookrightarrow X$ one can define a  ring homomorphism 
\[
	j_!:\,K(Y)\ \to\ K(X),
\] 
cf. \cite{AtSinger1} called  {\em push-forward map}.

The following special case of the {\em Bott periodicity} theorem, \cite{Atiyah68BottPeriodicity}, plays a crucial role in the definition of the topological index.

\begin{Thm}\label{T:Bottperiodicity} \textbf{\em (Bott periodicity)}
Let $i:\{pt\} \hookrightarrow \CC^N$ be an embedding of a point into $\CC^N$. Then the push-forward map
\[
	i_!:\,K(pt)\ \longrightarrow \ K(\CC^N)
\]
is an isomorphism. 
\end{Thm}

\subsection{The topological index}\label{SS:topindex}
Consider an embedding $M\hookrightarrow \RR^N$ (such an embedding always exists for large enough $N$, by the
Whitney embedding theorem). This embedding induces an embedding of the contangent
bundle of $M$
\[
    j:\, T^*M\ \hookrightarrow \RR^N\oplus \RR^N \ \simeq \ \CC^N.
\]
We will denote by $i:\{pt\}\hookrightarrow \CC^N$ an embedding of
a point into $\CC^N$. Using the push-forward in $K$-theory
introduced above we obtain the diagram
\begin{center}
\setlength{\unitlength}{.2in}
\begin{picture}(8,6)
\put(-2,4){$K(T^* M)$} \put(5,4){$K(\mathbb{C}^n)$}
\put(4.3,0){$K(\{pt\})$} \put(1.5,4.3){\vector(1,0){3.2}}
 \put(5.5,1.2){\vector(0,1){2.5}}
\put(3,5){$j_!$} \put(6,2){$i_!$}
\end{picture}
\end{center}

By Theorem~\ref{T:Bottperiodicity} we can invert it
and define the {\em topological index} as the map
\begin{equation}\label{E:tindex}
    \textrm{t-Ind} \ =\ i_!^{-1} j_! :K(T^* M) \rightarrow K(\{pt\})\ \simeq \ \mathbb{Z}.
\end{equation}

It is relatively easy to check that this map is independent of the choice of the number $N$ and the embedding $j:M\hookrightarrow \RR^N$.

We now introduce a group action into the picture. 
\subsection{Equivariant vector bundles}\label{SS:equivindexth}
Suppose a compact group $G$ acts on $M$. That means that to each $g\in G$ is assigned a diffeomorphism $\phi(g):M\to M$ 
such that $\phi(g_1)\circ\phi(g_2)= \phi(g_1g_2)$. 

\begin{Def}\label{D:eqivariantbundle}
A {\em $G$-equivariant vector bundle}  over $M$ is a vector bundle $\pi:E\to M$ together with an action 
\[
	g \ \mapsto \ \psi(g):\,E\ \to \ E, \qquad g\in G,
\]
of $G$ on $E$ such that 
\[
	\pi\circ\psi(g) \ =\ \phi\circ\pi(g).
\]
\end{Def}

When it does not lead to a confusion we often write $g\cdot{x}$ and $g\cdot{e}$ for $\phi(g)(x)$ and $\psi(g)(e)$ respectively (here $m\in M$ and $g\in G$). 

\begin{Def}\label{D:left regular}
Let  $E$ be a $G$-equivariant vector bundle over $M$. We define the action of  $G$ on the space $C^\infty(M,E)$ of smooth section of $E$ by the formula 
\begin{equation}\label{E:action on sections}
	g\ \mapsto \ l_g:\, C^\infty(M,E)\ \longrightarrow \ C^\infty(M,E), \qquad l_g(f)(x)\ := \ g\cdot f (g^{-1}\cdot x).
\end{equation}
Similarly, we define the action $l_g:L^2(M,E)\to L^2(M,E)$ on the space of square integrable sections of $E$.  
In this way $L^2(G)$ becomes a representation of $G$, called the {\em left regular representation}.

Suppose $E$ and $F$ are $G$-equivariant vector bundles over $M$. A differential operator 
\[
	\calD:\,C^\infty(M,E)\ \longrightarrow \ C^\infty(M,E),
\]
is called {\em $G$-invariant} if $g\cdot\calD= \calD\cdot{}g$. 
\end{Def}
\subsection{Equivariant K-theory}\label{SS:equivK}
If a compact group $G$ acts on $M$ one can define the {\em equivariant $K$-theory $K_G(X)$} as the set of the equivalence classes of $G$-equivariant maps between $G$-equivariant vector bundles. In particular, $G$-equivariant $K$-theory of a point is given by equivalence classes of pairs of finite dimensional representations of $G$. It should not be a surprise that
\begin{equation}\label{E:Kgpt}
    K_G(\{pt\}) \ = \ R(G),
\end{equation}
where $R(G)$ is the ring of characters of $G$, cf. Definition~\ref{D:R(G)}. The formula \eqref{E:tindex} generalizes easily to define a
$G$-equivariant topological index
\[
    \tInd_G:\, K_G(M)\ \to \ K_G(\{pt\})\ \simeq \ R(G).
\]
All the constructions introduced above readily generalize to the equivariant setting. 

\subsection{The Atiyah-Singer index theorem}\label{SS:ASindex}
Suppose a compact group $G$ acts on a {\em compact} manifold $M$. Let $E$ and $F$ be equivariant vector bundles over $M$ and let 
\[
	\calD:\,C^\infty(M,E)\ \to\ C^\infty(M,F)
\]
be a $G$-invariant elliptic operator. Then the leading symbol defines a $G$-equivariant map
\[
	\sig_L(\calD):\, \pi^*E\ \longrightarrow\ \pi^*F,
\]
which is invertible outside of the compact set $M\subset T^*M$. Thus $\sig_L\in K_G(T^*M)$ and its topological index $\tInd(\sig_L(\calD))$ is defined. The index theorem of Atiyah and Singer is the following result:

\begin{Thm}\label{T:eqindex} \textbf{\em (Atiyah-Singer)} 
$\Ind_G (\mathcal{D}) \ =\  \tInd_G(\sigma_L (\mathcal{D}))$.
\end{Thm}

\subsection{The case of an open manifold}\label{SS:openmanifold}
Suppose now that the manifold $M$ is not compact. The topological index map \eqref{E:tindex} is still defined. But the symbol of an elliptic operator $\calD$ does not define an element of $K(T^*M)$ since it is not invertible on a non-compact set $M\subset T^*M$.  Thus one can pose the following natural question

\subsection*{Question  1}
{\em Assume that a compact group $G$ acts on an {\em open} manifold $M$. Let $E$ and $F$ be $G$-equivariant vector bundles over $M$ and let $\sig:\pi^*E\to \pi^*F$ be an element of $K_G(T^*M)$. Find a $G$-invariant Fredholm differential operator $\calD:C^\infty(M,E)\to C^\infty(M,F)$ such that 
\[
	\Ind_G (\calD) \ =\  \tInd_G(\sig).
\]}

\medskip
To the best of our knowledge the answer to this question is unknown. However, in Section~\ref{S:indexnc} we present a partial answer to a certain generalization of this question.  For this we will need a generalization of the index theorem to so called {\em transversally elliptic} operators, which we will now discuss. 

\section{Transversal elliptic operators}\label{S:trelliptic}

In this section we discuss a generalization of the index theorem \ref{T:eqindex} to transversally elliptic operators due to Atiyah and Singer \cite{AtiyahTE}.

\subsection{A motivating example}\label{SS:MxG}
Suppose $N$ is a closed manifold and  let $\calD:C^\infty(N) \to C^\infty(N)$  be an elliptic operator. Let $G$ be a compact Lie group. Consider the manifold
\[
	M\ :=\ N \times G
\]
and let $G$ act on $M$ by 
\[
	g_1\cdot(x,g_2) \ := \ (x,g_1g_2), 	\qquad x\in N,\ g_1,g_2\in G.
\]

There is a natural extension of the differential operator $\calD$ to an operator
\[
	\tilde{\calD}:\,C^\infty(M) \ \longrightarrow \ C^\infty(M)
\]
defined as follows: Suppose that in a coordinate chart $(x_1\nek x_n)$ on $N$ the operator $\calD$ has the form 
\[
	\calD\ =\ \sum_{|\alpha|\leq k}a_\alpha(x)\,
	\frac{\partial^{\alp_1}}{\partial x_1^{\alp_1}}\cdots \frac{\partial^{\alp_n}}{\partial x_n^{\alp_n}}.
\]
Let $(y_1,\cdots, y_r)$ be a coordinate system on $G$. Then,  $(x_1\nek x_n,y_1\nek y_r)$ is a coordinate system on $M$. By definition the operator $\tilde{\calD}$ takes the form 
\[
	\tilde{\calD}\ :=\ \sum_{|\alp|\leq k}a_\alp(x)\, 
	\frac{\partial^{\alp_1}}{\partial x_1^{\alp_1}}\cdots \frac{\partial^{\alp_n}}{\partial x_n^{\alp_n}}.
\]

\begin{exercise}
Check that  $\tilde{\calD}$ is $G$-equivariant but not elliptic.
\end{exercise}

Notice that 
\begin{equation}\label{E:KertildeD}
	\Ker(\tilde\calD) \ = \ \Ker(\calD)\otimes L^2(G), \quad \Coker(\tilde\calD) \ = \ \Coker(\calD)\otimes L^2(G).
\end{equation}
Thus despite the fact that $\Ker(\tilde\calD)$ and $\Coker(\tilde\calD)$ are infinite dimensional they are sort of  ``manageable" and the index of $\tilde\calD$ can be defined as follows: 

By Peter-Weyl theorem \cite{Knapp01book} the left regular representation (cf. Definition~\ref{D:left regular})  $L^2(G)$ decomposes into direct sum of irreducible representations as 
\begin{equation}\label{E:Peter-Weyl}
	L^2(G)\ = \ \bigoplus_{V\in \Irr(G)}\,(\dim V)\,V.
\end{equation}
In particular every irreducible representation enters \eqref{E:Peter-Weyl} with finite multiplicity $\dim V$. 

Equations  \eqref{E:KertildeD} and \eqref{E:Peter-Weyl} suggest the index of $\tilde\calD$ can be defined as 
\begin{equation}\label{E:indexMxG}
	\Ind_G(\tilde\calD)\ = \ \bigoplus_{V\in \Irr(G)}\,(\dim V)\big(\dim\Ker(\calD)-\dim\Ker(\calD)\big)\,V.
\end{equation}
Notice that in contrast with \eqref{E:eqindex} the sum in the right hand side of \eqref{E:indexMxG} is infinite. Thus $\Ind(\tilde\calD)$ lies in the completed ring of characters $\hatR(G)$, cf. Remark~\ref{R:hatR}.

Roughly speaking the reason we are able to define a version of an index of $\tilde\calD$ is that though this operator is not elliptic on $M$ it is elliptic on the quotient $N=M/G$. The purpose of this section is to define the analogue of this situation when $M$ is not a product and the action of $G$ on $M$ is not free. 

\subsection{The transversal cotangent bundle}\label{SS:TGM}
Suppose a compact group $G$ acts on a smooth manifold $M$. Recall that an {\em orbit} of a point $x\in M$ is the set 
\[
	\calO(x)\ := \ \big\{\,g\cdot x\big|\, g\in G\,\big\}.
\]
This is a smooth submanifold of $M$. 

We say that a cotangent vector $\xi\in T^*M$ is {\em perpendicular to the orbits of $G$} if for $x\in M$ and any tangent vector $v\in T\calO(x)$ we have $\xi(v)=0$. 

\begin{Def}\label{D:TGM}
The {\em transversal cotangent bundle $T^*_GM$} is defined by
\begin{equation}
	T_G^*M\ :=\ \big\{\xi \in T^*M\big|\, \xi \ \text{\em  is perpendicular to the orbits of } G\,\big\}.
\end{equation}
We set $T^*_{G,x}M:= T^*_GM\cap T^*_xM$.
\end{Def}

\subsection{The analytical index of transversally elliptic operators}\label{SS:trelliptic}
We now introduce the class of operators which generalizes the example considered in  Section~\ref{SS:MxG}.
\begin{Def}\label{D:trelliptic}
Suppose $E$ and $F$ are $G$-equivariant vector bundles over $M$. A $G$-invariant differential operator
\[
	\calD:\,C^\infty(M,E)\ \longrightarrow\ C^\infty(M,F)
\]
is called {\em transversally elliptic} if its leading symbol $\sig_L(\calD)(x,\xi)$ is invertible for every non-zero $0\not=\xi \in T^*_GM$. 
\end{Def}
Notice that every elliptic operator is transversally elliptic. Also the operator $\tilde\calD$ from Section~\ref{SS:MxG} is transversally elliptic. 

\begin{Thm}\label{T:trelliptic} \textbf{\em (Atiyah-Singer \cite[Lemma~2]{AtiyahTE})}\ 
Suppose $\calD$ is a transversally elliptic operator on a {\em compact} manifold $M$. Then, as the representations of $G$, the spaces $\Ker(\calD)$ and $\Coker(\calD)$ can be decomposed into a direct sum of irreducible representation in which every irreducible representation appears finitely many times:
\begin{equation}\label{E:trellipticKer}
	\Ker(\calD)\ = \ \bigoplus_{V\in \Irr(G)} m_V^+V, \qquad \Coker(\calD)\ = \ \bigoplus_{V\in \Irr(G)} m_V^-V.
\end{equation}
\end{Thm}
Notice that the sums in \eqref{E:trellipticKer} are, in general, infinite. But the numbers $m_V^\pm\in \ZZ_{\ge0}$ are finite. 
\begin{Def}\label{D:trindex}
Let $\calD$ be a  transversally elliptic operator on a {\em compact} manifold $M$. The (analytical) index of $\calD$ is defined as 
\begin{equation}\label{E:trindex}
	\Ind_G(\calD)\ := \ \bigoplus_{V\in \Irr(G)} (m_V^+-m_V^-)\,V\ \in \ \hatR(G),
\end{equation}
where the numbers $m_V^\pm\in \ZZ_{\ge0}$ are defined in \eqref{E:trellipticKer} and $\hatR(G)$ stands for the completed ring of characters, cf. Remark~\ref{R:hatR}.
\end{Def}

The index \eqref{E:trindex} possesses many properties of the index of elliptic operators. In particular and analogue of Theorem~\ref{T:homotinvindex} holds.

\subsection{The transversal K-theory and the topological index}\label{SS:trKtheory}
Notice that in general $T_G^*M$ is not a manifold since the dimension of the fibers $T^*_{G,x}M$ might depend on $x$. But it is a topological space and one can define the K-theory $K_G(T^*_GM)$ as the set of equivalence classes of pairs of vector bundles over $T^*_GM$ exactly the same way as we did in Section~\ref{SS:equivK}. 

The topological index 
\begin{equation}\label{E:trtindex}
	\tInd_G:\,K_G(T^*_GM)\ \longrightarrow\ \hatR(G)
\end{equation}	
is still defined, but  the image lies not in the ring of characters and the definition is more involved. In fact, for compact manifold $M$ if $\sig\in K_G(T^*_GM)$ one just chooses a $G$-invariant operator $\calD$  with $\sig_L(D)=\sig$ and sets
\begin{equation}\label{E:AStelliptic}
	\tInd_G(\sig)\ := \ \Ind_G(\calD).
\end{equation}
(see \cite{BerlineVergne96} for a more topological construction). For non-compact $M$ the topological index is defined in  \cite[\S3]{Paradan01}. As in Section~\ref{SS:openmanifold}, if $M$ is not compact the symbol of a transversally elliptic operator does not define an element of $K_G(T^*_GM)$. The index of such an operator is not defined since, in general, the irreducible representations of $G$ appear in the kernel and the cokernel of $\calD$ with infinite multiplicities. So as in  Section~\ref{SS:openmanifold} a natural question arises

\subsection*{Question  2}
{\em Assume that a compact group $G$ acts on an {\em open} manifold $M$. Let $E$ and $F$ be $G$-equivariant vector bundles over $M$ and let $\sig:\pi^*E\to \pi^*F$ be an element of $K_G(T^*_GM)$. Find a $G$-invariant  differential operator $\calD:C^\infty(M,E)\to C^\infty(M,F)$ such that each irreducible representations of $G$ appear in the kernel and the cokernel of $\calD$ with finite multiplicities and
\[
	\Ind_G (\calD) \ =\  \tInd_G(\sig)\ \in \hatR(G).
\]}
\smallskip
Of course, this question is even harder than the corresponding question for usual K-theory (cf. Question~1 in Section~\ref{SS:openmanifold}) and the answer to this question is unknown. However in Section~\ref{S:indexnc} we will present an answer for Question~2 in a special case. To present this answer we need to introduce the notions of Clifford action and Dirac-type operators, which we do in the next section.

\section{Dirac-type operators}\label{S:Dirac}

In this section we define the notions of Clifford bundle and generalized Dirac operator. 
\subsection{Clifford action}\label{SS:Clifford action}
Let $V$ be a finite dimensional vector space over $\RR$ endowed with a scalar product $\<\cdot,\cdot\>$. 

\begin{Def}\label{D:Clifford action}
A {\em Clifford action} of $V$ on a complex vector space $W$ is a linear map 
\[
	c:\,V\ \to\ \End(W)
\] 
such that for any $v\in V$ we have
\begin{equation}\label{E:Clifford action}
	c(v)^2\ = \ -|v|^2\,\Id.
\end{equation}
\end{Def}

\begin{exercise}
Show that a linear map  $c:V\to \End(W)$  is a Clifford action if and only if 
\[
	c(v)\,c(u)+c(u)\,c(v) \ = \ -2\<v,u\>\,\Id, \qquad\text{for all}\ u,v\in V.
\]
\end{exercise}

\begin{exercise} \label{Exer:Pauli}
Consider the Pauli matrices:
\begin{equation}\label{E:Pauli}
	\sig_1\ = \  \begin{pmatrix} 0&1\\1&0\end{pmatrix}, \qquad
	\sig_2\ = \  \begin{pmatrix} 0&-i\\i&0\end{pmatrix}, \qquad
	\sig_3\ = \  \begin{pmatrix} 1&0\\0&-1\end{pmatrix}.
\end{equation}
Show that the map
\begin{equation}\label{E:spinoraction}
	c:\,\RR^3 \ \longrightarrow \ \End(\CC^2), \quad c(x_1,x_2,x_3) \ := \ \frac1i\,\big(\,x_1\sig_1+x_2\sig_2+x_2\sig_3\,\big)
\end{equation}
defines a Clifford action of $\RR^3$ on $\CC^2$. 
\end{exercise}

\begin{exercise}\label{Exer:exterior}
Let $V$ be a real vector space endowed with a scalar product and let $V^\CC= V\otimes_\RR\CC$ be its complexification. Let
\[
	W\ = \ \Lam^*(V^\CC)
\]
be the exterior algebra of $V^\CC$. For $v\in V$  we denote $\eps(v):W\to W$ the exterior multiplication by $v$
\[
	\eps(v)\,\alp\ := \ v\wedge\alp, \qquad \alp\in W.
\]
Also using the scalar product on $V$ we identify $v$ with an element of the dual space $V^*$ and denote by $\iot_v:W\to W$ the interior multiplication by $v$. Set 
\begin{equation}\label{E:DeRham Clifford}
	c(v)\ :=\ \eps_v\ - \ \iot_v:\, W\ \longrightarrow\ W.
\end{equation}
Show that the map $c:v\mapsto c(v)$ defines a Clifford action of $V$ on $W$. 
\end{exercise}

\subsection{A Clifford bundles}\label{SS:Clifford bundle}
Let now $M$ be a smooth manifold endowed with a Riemannian metric $g$. Then for every $x\in M$ the tangent space $T_xM$ and the cotangent space $T^*_xM$ are endowed with a scalar product. 

\begin{Def}\label{D:Clifford bundle}
A {\em Clifford bundle} over $M$ is a complex vector bundle $E\to M$ over $M$ together with a bundle map
\begin{equation}\label{E:Clifford bundle}
	c:\,T^*M \ \longrightarrow \  \End(E),
\end{equation}
such that for any $v\in T^*M$ we have $c(v)^2=-|v|^2\Id$.
\end{Def}
In other words we assume that for every $x\in M$ there is given a Clifford action 
\[
	c:\,T_x^*M\ \to\ \End(E_x)
\]
of the cotangent space $T_x^*M$ on the fiber $E_x$ of $E$ and that this action depends smoothly on $x$.

\begin{exam}\label{Exam:clDirac}
Set $M=\RR^3$ and $E=\RR^3\times\CC^2$. Then the projection on the first factor makes $E$ a vector bundle over $M$. The action  \eqref{E:spinoraction} defines a structure of a Clifford bundle on $E$.
\end{exam}

\begin{exam}\label{Exam:DeRham Clifford}
Let $M$ be a Riemannian manifold and let $E=\Lam^*(T^*M\otimes\CC)$. Note that the space of sections of $E$ is just the space $\Ome^*(M)$ of complex valued differential forms on $M$. The formula \eqref{E:DeRham Clifford} defines a Clifford bundle structure on $E$, such that for $v\in T^*M$ and $\ome\in \Ome^*(M)$ we have
\[
	c(v)\,\ome\ = \ v\wedge \ome\ - \ \iot_v\,\ome.
\]
\end{exam}


\subsection{A Clifford connection}\label{SS:Clifford connection}
The Riemannian metric on $M$ defines a connection on $T^*M$ called the {\em Levi-Civita connection} and denoted by $\n^{LC}$.

\begin{Def}\label{D:Clifford connection}
A connection $\n:C^\infty(M,E)\to \Ome^1M,E)$  on a Clifford bundle $E$ is called a {\em Clifford connection} if 
\begin{equation}
	\n_u\big(c(v)s\big)\ = \ c(\n^{LC}_uv)\,s\ + \ c(v)\,\n_us, 
\end{equation}
for all $u\in TM, \ v\in T^*M,\ s\in C^\infty(M,E)$.
\end{Def}

\subsection{A generalized Dirac operator}\label{SS:Dirac}
We are now ready to define a notion of (generalized) Dirac operator. 

\begin{Def}\label{D:Dirac}
Let $E$ be a Clifford bundle over a Riemannian manifold $M$ and let $\n$ be a Clifford connection on $E$. The {\em (generalized) Dirac operator} is defined by the formula
\begin{equation}\label{E:Dirac}
	\calD\ := \ \sum_{j=1}^n\,c(e_i)\,\n_{e_i}:\, C^\infty(M,E) \ \longrightarrow \ C^\infty(M,E),
\end{equation}
where $e_1,\ldots,e_n$ is an orthonormal basis of $TM$.%
\footnote{As above, we use the Riemannian metric to identify $TM$ and $T^*M$.}
\end{Def}

\begin{exercise}
Show that the operator \eqref{E:Dirac} is independent of the choice of the orthonormal basis $e_1,\ldots,e_n$.
\end{exercise}
\begin{exercise}\label{Exer:Diracelliptic}
a. \ Show that the leading symbol generalized Dirac operator is given by 
\begin{equation}\label{E:Diracsymbol}
	\sig_L(\calD)(x,\xi) \ = \ i\,c(\xi).
\end{equation}

b. \ Prove that the generalized Dirac operator is elliptic. 
\end{exercise}

\begin{exam}
In the situation of Example~\ref{Exam:clDirac} let $\n$ be the standard connection, i.e.,
\[
	\n_{\partial/\partial x_i}\ = \ \frac{\partial}{\partial x_i}.
\]
Then \eqref{E:Dirac} is equal to the classical Dirac operator
\begin{equation}\label{E:clDir}
	\calD \ = \  \begin{pmatrix} 0&1\\1&0\end{pmatrix}\, \frac1i\frac{\partial}{\partial x_1} \ + \
	 \begin{pmatrix} 0&-i\\i&0\end{pmatrix}\,  \frac1i\frac{\partial}{\partial x_2} \ + \  
	 \begin{pmatrix} 1&0\\0&-1\end{pmatrix}\,  \frac1i\frac{\partial}{\partial x_3}.
\end{equation}
\end{exam}
\begin{exercise}
Compute the square of the operator \eqref{E:clDir}.
\end{exercise}

\begin{exam}
Let $M$ be a Riemannian manifold and let $d:\Ome^*(M,E)\to \Ome^{*+1}(M,E)$ denote the de Rham differential. The Riemannian metric induces a scalar product $\<\cdot,\cdot\>$ on the space $\Ome^*(M,E)$ of differential forms. Let 
\[
	d^*:\,\Ome^*(M,E)\ \to\  \Ome^{*-1}(M,E)
\] 
be the adjoint on $d$ with respect to this scalar product. Thus we have
\[
	\<d\alp,\bet\>\ =\ \<\alp,d^*\bet\>, \qquad\text{for all}\quad \alp,\bet\in \Ome^*(M,E).
\]

Let $E=\Lam^*(T^*M\otimes\CC)$. Then the space of sections $C^\infty(M,E)= \Ome^*(M,E)$. The Riemannian metric defines a canonical connection on $M$, called the Levi-Civita connection, cf. \cite{Taubes-DifGeom}.  We denote this connection by $\n^{LC}$. One can show, cf. \cite[Proposition~3.53]{BeGeVe}, that this is a Clifford connection and that the corresponding Dirac operator \eqref{E:Dirac} is given by
\begin{equation}\label{E:deRhamDirac}
	\calD\ = \ d\ + \ d^*.
\end{equation}
If $M$ is a compact manifold then, \cite{Warner}, the kernel of $\calD$ is naturally isomorphic to the de Rham cohomology on $M$:
\begin{equation}\label{E:deRhamcoh}
	\Ker(\calD)\ \simeq \ H^*(M).
\end{equation}
\end{exam}

\subsection{A grading}
In the index theory one deals with graded Clifford bundles. 
\begin{Def}\label{D:grading}
A {\em grading} on a Clifford bundle $E$ is a decomposition 
\begin{equation}\label{E:grading}
	E\ =\ E^+\oplus E^-,
\end{equation}
such that for every $v\in T^*M$ we have
\[
	c(v):\,E^\pm\ \longrightarrow\ E^\mp.
\]
A Clifford bundle with a grading is called a {\em graded Clifford bundle}.

We say that a Clifford connection on  a graded Clifford bundle {\em preserves the grading} if  for each  for each $v\in T^*M$, 
\[
	\n_v:\,C^\infty(M,E^\pm) \ \longrightarrow \ C^\infty(M,E^\pm).
\]
\end{Def}
When speaking about a Clifford  connection $\n$ on a graded Clifford bundle,  we will always assume that it preserves the grading. 

\begin{exercise}\label{Exer:grading}
Let $E=E^+\oplus E^-$ be a graded Clifford bundle and let $\n$ be a Clifford connection on $E$ which preserves the grading. Show that the corresponding Dirac operator satisfies 
\[
	\calD :\,C^\infty(M,E^\pm) \ \longrightarrow \ C^\infty(M,E^\mp).
\]
\end{exercise}

We denote the restriction of $\calD$ to $C^\infty(M,E^+)$ (respectively $C^\infty(M,E^-)$) by $\calD^+$ (respectively $\calD^-$).
It follows from the above exercise that with respect to the splitting \eqref{E:grading} the Dirac operator $\calD$ can be written as 
\begin{equation}\label{E:gradedDirac}
	\calD\ = \ \begin{pmatrix}0&\calD^-\\\calD^+&0\end{pmatrix}.
\end{equation}
To save the space we often write \eqref{E:gradedDirac} as $\calD=\calD^+\oplus\calD^-$. 
When a Dirac operator is presented in this form we refer to it as a {\em graded Dirac operator}.  In the index theory one usually considers the index of the operator $\calD^+:C^\infty(M,E^+)\to C^\infty(M,E^-)$.

\newcommand{\even}{\operatorname{even}}\newcommand{\odd}{\operatorname{odd}}
\begin{exam}
Let $M$ be a Riemannian manifold and let $E=\Lam^*(T^*M\otimes\CC)$. Set 
\[
	E^+\ := \ \bigoplus_{j \even}\Lam^j(T^*M\otimes\CC), \qquad E^-\ := \ \bigoplus_{j \odd}\Lam^j(T^*M\otimes\CC).
\]
Then $E=E^+\oplus E^-$, i.e., we obtain a grading on $E$. The Levi-Civita connection preserves this grading. Thus \eqref{E:deRhamDirac} becomes a graded operator. This graded Dirac operator is called the {\em de Rham-Dirac operator}. 
\end{exam}
\begin{exercise}
Suppose that the manifold $M$ is compact.  Show that the index $\Ind(\calD^+)$ of the de Rham-Dirac operator is equal to the Euler characteristic of $M$:
\[
	\Ind(\calD^+) \ = \ \sum_{j=0}^n\,(-1)^n\dim H^j(M).
\]
\end{exercise}

\begin{exam}
Suppose the dimension of $M$ is even, $\dim M= n=2l$. There is another natural grading on $E=\Lam^*(T^*M\otimes\CC)$ defined as follows. Let $*:E\to E$ denote the Hodge star operator \cite[\S19.1]{Taubes-DifGeom}. Define the {\em chirality operator} $\Gam:E\to E$ by 
\[
	\Gam\ome\ := \ i^{p(p-1)+l}*\ome, \qquad \ome\in \Lam^p(M)\otimes\CC.
\]
One can show, \cite[\S19.1]{Taubes-DifGeom}, that $\Gam^2=1$. It follows that the spectrum of the operator $\Gam$ is the set $\{1,-1\}$. Let $E^+$ (respectively, $E^-$) denote the eigenspace of $\Gam$ corresponding to the eigenvalue +1 (respectively, -1). The sections of $E^+$ (respectively $E^-$) are called the {\em self-dual} (respectively, {\em anti-selfdual}) differential forms. Then $E=E^+\oplus E^-$ and the Levi-Civita connection preserves this grading. Hence, \eqref{E:deRhamDirac} becomes a graded operator. This graded Dirac operator is called the {\em signature operator} and the index of $\calD^+$ is called the {\em signature} of $M$. The signature of $M$ can also be computed in topological terms. The study of the index of the signature operator by Friedrich Hirzebruch \cite{Hirzebruch66book} was one of the main motivation for the  Atiyah-Singer work on index theory. 
\end{exam}

\subsection{The group action}\label{SS:groupaction Dirac}
Let $E=E^+\oplus E^-$ be a graded Clifford bundle over a Riemannian manifold $M$. Let $\n$ be a Clifford connection which preserves the grading and let $\calD=\calD^+\oplus\calD^-$ be the corresponding graded Dirac operator. Suppose that a compact group $G$ acts on $M$ and $E$, preserving the Riemannian metric, the connection $\n$, and the grading on $E$. Then the operators $\calD^\pm$ are $G$-equivariant. In particular, if the manifold $M$ is compact, we can consider the equivariant index
\[
	\Ind_G(\calD^+)\ \in \ R(G).
\]

\section{Index theory on open $G$-manifolds}\label{S:indexnc}

We are now ready to define the class of transversally elliptic symbols for which we are able to answer Question~2 of Section~\ref{SS:trKtheory}. The section is based on the results of  \cite{Br-index}.

\subsection{The settings}\label{SS:settingsnc}
Recall that a Riemannian manifold $M$ is called {\em complete} if it is complete as a metric space. Suppose $M$ is a complete Riemannian manifold, on which a compact Lie group $G$ acts by isometries. To construct our index theory of certain generalized Dirac operators on $M$ we need an additional structure on $M$, namely a $G$-equivariant map 
\begin{equation}\label{E:map v}
	\bfv:\,M\ \to\ \grg\ =\ \Lie G.
\end{equation} 
Such  a map induces a vector field $v$ on $M$ defined by the formula
\begin{equation}\label{E:v}
    v(x) \ := \ \frac{d}{dt}\Big|_{t=0}\, \exp{(t\bfv(x))}\cdot x.
\end{equation}

\begin{exam}\label{Exam:generatingvf}
Let 
\[
	G\ =\ S^1\ =\ \big\{\,z\in\CC:\, |z|=1\,\big\}
\]
be the circle group. Then the Lie algebra $\grg=\Lie{}S^1$ is naturally isomorphic to $\RR$. Suppose $G$ acts on a complete Riemannian manifold $M$. Let $\bfv:M\to\grg$ be the constant map $\bfv(x)\equiv1$. The corresponding vector field $v$ is called the {\em generating vector field} for the action of $G=S^1$, since it completely determines the action of $G$. 
\end{exam}
\begin{exercise}
Let $G=S^1$ acts on $M=\CC$ by multiplications: $(e^{it},z)\mapsto e^{it}z$. Compute the generating vector field $v(x)$.
\end{exercise}

\subsection{Tamed $G$-manifolds}\label{SS:tamed}
Throughout this section we will make the following 
\subsection*{Assumption} There exists a compact subset $K\subset M$, such that 
\begin{equation}\label{E:assumption}
	v(x)\ \not=\ 0, \qquad\text{for all}\quad x\not\in K.
\end{equation}

\begin{Def}\label{D:tamed}
A map \eqref{E:map v} satisfying \eqref{E:v}  is called a {\em taming map}. The pair $(M,\bfv)$, where $M$ is a complete Riemannian manifold and $\bfv$ is a taming map,  is  called a {\em tamed $G$-manifold}.
\end{Def}

\subsection{An element of $K_G(T^*_GM)$}\label{SS:elementKG}
Suppose now that $(M,\bfv)$ is a tamed $G$-manifold and let $E=E^+\oplus{}E^-$ be a $G$-equivariant graded Clifford bundle over $M$. Let $T^*_GM$ denote the transversal cotangent bundle to $M$, cf. Definition~\ref{D:TGM}, and let $\pi:T^*_GM\to M$ denote the projection. Consider the pull-back bundle $\pi^*E=\pi^*E^+\oplus\pi^*E^-$ over $T^*_GM$. 

Using the Riemannian metric on $M$ we can identify the tangent and cotangent vectors to $M$. Thus we can consider the vector $v(x)$ defined in \eqref{E:assumption} as an element of $T^*M$. Then for $x\in M$, $\xi\in T^*M$ we can consider the map
\begin{equation}\label{E:sigxi-v}
	c\big(\xi+v(x)\big):\, \pi^*E^+_{x,\xi}\ \to\  \pi^*E^-_{x,\xi}.
\end{equation}
The collection of those maps for all $(x,\xi)\in T^*_GM$ defines a bundle map 
\[
	c(\xi+v):\,\pi^*E^+\ \to\ \pi^*E^-.
\]

\begin{exercise}\label{Exer:sigxi-v}
Show that the map $c\big(\xi+v(x)\big)$ is invertible for all $(x,\xi)\in T^*_GM$ such that $v(x)\not=0$, $\xi\not=0$. Conclude that the bundle map \eqref{E:sigxi-v} defines an element of the K-theory $K_G(T^*_GM)$.
\end{exercise}

Hence, we can consider the topological index 
\begin{equation*}
	\tInd_G\big(c(\xi+v)\big) \ \in \  \hatR(G).
\end{equation*}

\begin{Def}
Suppose $E=E^+\oplus E^-$ is a $G$-equivariant Clifford bundle over a tamed $G$-manifold $(M,\bfv)$. We refer to the pair $(E,\bfv)$ as a {\em tamed Clifford bundle over $M$}. The topological index $\tInd_G(E,\bfv)$ of a tamed Clifford bundle is defined by 
\begin{equation}\label{E:indexEv}
	\tInd_G(E,\bfv)\ := \ \tInd_G\big(c(\xi+v)\big).
\end{equation}
\end{Def}

This index was extensively studied by M.~Vergne \cite{Vergne96} and P.-E.~Paradan \cite{Paradan01,Paradan03}.
 Our purpose is to construct a Fredholm operator, whose analytical index is equal to $\tInd_G(E,\bfv)$.

\subsection{The Dirac operator}\label{Diracv}
Suppose $\n$ is a $G$-equivariant Clifford connection on $E$ which preserves the grading and let $\calD=\calD^+\oplus\calD^-$ be the corresponding Dirac operator.  By Exercise~\ref{Exer:Diracelliptic} the operator $\calD$ is elliptic and its leading symbol is equal to $c(\xi)$. However, if $M$ is not compact, the operator $\calD$ does not have to be Fredholm and its index is not defined.

\subsection{A rescaling of $v$} \label{SS:rescaling}
Our definition of the index uses certain rescaling of the vector field  $v$. By this we mean the product $f(x)v(x)$,
where $f:M\to[0,\infty)$ is a smooth positive $G$-invariant function. Roughly speaking, we demand that $f(x)v(x)$ tends to
infinity ``fast enough" when $x$ tends to infinity. The precise conditions we impose on $f$ are quite technical,
cf. Definition~2.6 of \cite{Br-index}, and depends on the geometry of $M$, $E$ and $\n$.  If a function $f$ satisfies these conditions we call it {\em admissible for the quadruple $(M,E,\n,\bfv)$}.  The index which we are about to define  turns out to be independent of the concrete choice of $f$. It is important, however, to know that at least one admissible function exists. This is proven in Lemma~2.7 of \cite{Br-index}.

\subsection{The analytic index on non-compact manifolds}\label{SS:noncomind}
Let $f$ be an admissible function for $(M,E,\n,\bfv)$. Consider the {\em deformed Dirac operator}
\begin{equation}\label{E:Dv}
        \calD_{fv} \ = \ \calD \ + \ {i}c(fv),
\end{equation}
and let $\calD^+_{fv}$ denote the restriction of $\calD_{fv}$ to $C^\infty(M,E^+)$.  This operator is elliptic, but since the manifold $M$ is not compact it is not Fredholm. In fact, both, the kernel and the cokernel of $\calD^+_{fv}$, are infinite dimensional. However they have an important property, which we shall now describe. First, recall from Exercise~\ref{Exer:equivariant kernel} that, since the operator $\calD^+_{fv}$ is $G$-equivariant, the group $G$ acts on $\Ker(\calD^+_{fv})$ and $\Coker(\calD^+_{fv})$.

\begin{Thm}\label{T:finite}
Suppose $f$ is an admissible function.  Then the kernel and the cokernel of the deformed Dirac operator $D^+_{fv}$
decompose into an infinite direct sum
\begin{equation}\label{E:finite}
        \Ker \calD^+_{fv} \ = \ \sum_{V\in \Irr(G)}\, m^+_V\cdot V, \quad  \Coker \calD^+_{fv} \ = \ \sum_{V\in \Irr(G)}\, m^-_V\cdot V,
\end{equation}
where $m_V^\pm$ are non-negative integers. In other words, each irreducible representation of $G$ appears in
$\Ker \calD^\pm_{fv}$ with finite multiplicity.
\end{Thm}

This theorem allows us to define the index of the operator $\calD^+_{fv}$:
\begin{equation}\label{E:index Dfv}
	\Ind_G(\calD^+_{fv}) \ = \ \sum_{V\in \Irr(G)}\, (m^+_V-m_V^-)\cdot V \ \in \ \hatR(G).
\end{equation}
The next theorem states that this index is independent of all the choices. 
\begin{Thm}\label{T:indDfv}
For an irreducible representation $V\in \Irr(G)$ let $m_V^\pm$ be defined by \eqref{E:finite}. Then the differences $m^+_V-m^-_V$ $(V\in\Irr(G))$ are independent of the choices of the admissible function $f$ and the $G$-invariant Clifford connection on $E$, used in the definition of $\calD$.
\end{Thm}

\begin{Def}\label{D:indexDv}
We refer to the pair $(\calD,\bfv)$ as a {\em tamed Dirac operator}. The analytical index of a tamed Dirac operator is defined by the formula 
\begin{equation}\label{E:indexDv}
	\Ind_G(\calD,\bfv)\ := \ \Ind_G(\calD^+_{fv}),
\end{equation}
where $f$ is any admissible function for $(M,\bfv)$.
\end{Def}

\begin{exercise}
Let $(M,\bfv)$ be the tamed $G=S^1$-manifold defined in Example~\ref{Exam:generatingvf}. Let 
\[
	E^\pm\ =\ M\times\RR
\]
be two line bundles over $M$. Define the Clifford action of $T^*M$ on $E=E^+\oplus{}E^-$ by the formula
\[
	c(\xi)\ := \  \frac1i\,\big(\,\xi_1\,\sig_1+ \xi_2\,\sig_2\,\big), \qquad \xi=(\xi_1,\xi_2)\in T^*M\simeq \RR^2,
\]
where $\sig_1$ and $\sig_2$ are the first two Pauli matrices, cf. \eqref{E:Pauli}.

The $G$ action on $M$ lifts to $E^\pm$ so that 
\[
	e^{it}\cdot(z,\nu)\ := \ (e^{it}z,\nu), \qquad z\in M, \ t,\nu\in \RR.
\]
We endow the bundle $E$ with the trivial connection.

a. \ Show that the above construction defines a structure of a $G$-equivariant graded Clifford module on $E$. 

b. \ One can show that $f\equiv1$ is an admissible function for $E$. Compute the operators $\calD_v$ and $\calD_v^2$. 

c. \ Find the kernel and cokernel of $\calD_v$. 

d. \ Compute $\Ind_G(\calD,\bfv)$. 
\end{exercise}

\subsection{The index theorem}\label{SS:index theorem}
We are now ready to formulate the following analogue of the Atiyah-Singer index theorem for non-compact manifolds:

\begin{Thm}\label{T:ASnoncompact}
Suppose $E=E^+\oplus E^-$ is a graded $G$-equivariant Clifford bundle over a tamed $G$-manifold $(M,\bfv)$. Let $\n$ be a $G$-invariant Clifford connection on $E$ which preserves the grading and let $\calD$ be the corresponding Dirac operator. Then for any 
\begin{equation}\label{E:ASnoncompact}
	\Ind(\calD,\bfv)\ =\ \tInd(E,\bfv).
\end{equation}
\end{Thm}

\begin{exercise}\label{Exer:Dvcompact}
Show that if the manifold $M$ is compact, then
\begin{equation}\label{E:Dvcompact}
	\Ind_G(\calD,\bfv)\ = \ \Ind_G(\calD).
\end{equation}
\end{exercise}


\subsection{Properties of the index on a non-compact manifold}\label{SS:propertiesInd}
The index \eqref{E:indexDv} has many properties similar to the properties of the index of an elliptic operator on a compact manifold. It satisfies an analogue of the Atiyah-Segal-Singer fixed point theorem, Guillemin-Sternberg  "quantization commutes with reduction" property, etc. It is also invariant under a certain type of cobordism. The description of all these properties lies beyond the scope of these lectures and we refer the reader to \cite{Br-index} for details.  We will finish with mentioning just one property -- the gluing formula -- which illustrates how the non-compact index can be used in the study of the usual index on compact manifolds. 

\subsection{The gluing formula}\label{SS:gluing}
Let $(M,\bfv)$ be a tamed $G$-manifold. Suppose $\Sig\subset M$ is a smooth $G$-invariant hypersurface in $M$ such that $M\backslash{\Sig}$ is a disjoint union of two open manifolds $M_1$ and $M_2$:
\[
	M\backslash \Sig \ = \ M_1 \sqcup M_2.
\]
 For simplicity, we assume that $\Sig$ is compact. Assume also that the vector field $v$ induced by $\bfv$ does not vanish anywhere on $\Sig$. Choose a $G$-invariant complete Riemannian metric on $M_j$ ($j=1,2)$. Let $\bfv_j$ denote the restriction of $\bfv$ to $M_j$. Then $(M_j,\bfv_j)$ $(j=1,2)$ are tamed $G$-manifolds.

Suppose that $E=E^+\oplus{}E^-$ is a $G$-equivariant graded Clifford module over $M$. Denote by $E_j$ the
restriction of $E$ to $M_j$ ($j=1,2)$. Let $\calD_j$  $(j=1,2)$ denote the restriction of $\calD$ to $M_j$.

\begin{Thm}\label{T:gluing}
In the situation described above
\begin{equation}\label{E:gluing}
	\Ind_G(\calD,\bfv)\ = \ \Ind_G(\calD_1,\bfv_1) \ +\ \Ind_G(\calD_2,\bfv_2).
\end{equation}
\end{Thm}

In view of Exercise~\ref{Exer:Dvcompact}, one can use Theorem~\ref{T:gluing} for studying the index of an equivariant Dirac operator on a compact manifold. This is done by cutting a compact manifold $M$ along a $G$-invariant hypersurface $\Sig$ into two non-compact, but topologically simpler manifolds $M_1$ and $M_2$. We refer the reader to \cite{Br-index} for examples of different applications of this idea.

\def\cprime{$'$} \def\cprime{$'$} \newcommand{\noop}[1]{} \def\cprime{$'$}
\providecommand{\bysame}{\leavevmode\hbox to3em{\hrulefill}\thinspace}
\providecommand{\MR}{\relax\ifhmode\unskip\space\fi MR }
\providecommand{\MRhref}[2]{%
  \href{http://www.ams.org/mathscinet-getitem?mr=#1}{#2}
}
\providecommand{\href}[2]{#2}

\end{document}